\documentclass[11pt]{amsart} 
\usepackage{amssymb,amsmath,latexsym,enumerate,graphicx,bbm,mathptmx,microtype}
\usepackage{float}
\usepackage{hyperref}
\usepackage{tabularx}
\usepackage{tabularx}
\usepackage{array}

\usepackage[a4paper]{geometry}
\numberwithin{equation}{section}
\newtheorem{theorem}{Theorem}[section]

\newtheorem{proposition}{Proposition}[section]

\newtheorem{exam}{Example}

\newtheorem*{rem}{Remark}
\newenvironment{remark}{\begin{rem}\rm}{\end{rem}}

\newcommand\commentout[1]{}

\usepackage{tikz}
\usetikzlibrary{mindmap}

\begin{document}
\date{\today}
\title[\tiny{THE SIMPLEST QUARTIC FIELDS.}]{INDEX,  PRIME IDEAL FACTORIZATION  IN SIMPLEST QUARTIC FIELDS and counting their discriminants.}
\author{ 
 Mohammed Seddik}
\address{ 
 Mohammed Seddik\\
Universit\'e d'\'Evry Val d'Essonne\\
Laboratoire de MathŽmatiques et Mod\'elisation d'\'Evry (UMR 8071)\\
I.B.G.B.I., 23 Bd. de France, 91037 \'Evry Cedex, France\\}
 
  \email{mohammed.seddik@univ-evry.fr}
\thanks{\footnotesize{$^*$ \it Corresponding Author: Mohammed Seddik },\  mohammed.seddik@univ-evry.fr}

\begin{abstract} We consider the simplest quartic number fields $\mathbb{K}_m$ defined by the irreducible quartic polynomials
$$x^4-mx^3-6x^2+mx+1,$$ where $m$ runs over the positive rational integers such that the odd part of $m^2+16$ is squarefree. 
In this paper, we study the common index divisor $I(\mathbb K_m)$ and determine explicitly the prime ideal decomposition for any prime number in any  simplest quartic number fields $\mathbb{K}_m$. 
On the other hand, 
we establish an asymptotic formula for the number of simplest quartic fields with discriminant $\leq x$ and given index. 
\end{abstract}




\maketitle

\section{\textbf {Introduction and preliminaries}}

Let $\mathbb K$ be a number field.
We start with given most important known results on the common factor of indices $I(\mathbb K)$, the prime ideal factorization and the asymptotic formula for number fields with discriminant $\leq x$.
\subsection{Index $I(\mathbb K)$ for number fields $\mathbb K$}
Let $\mathbb{K}$ be a number field of degree $n$ over $\mathbb{Q}$ and let $\mathbb{O_K}$ be its ring of integers. Denote by $\widehat{\mathbb{O}}_\mathbb{K}$ the set of primitive elements of $\mathbb{O_K}$. For any $\theta\in \mathbb{O_K}$ we denote $F_\theta(x)$ the characteristic polynomial of $\theta$ over $\mathbb{Q}$. Let $D_\mathbb{K}$ be the discriminant of $\mathbb{K}$. It is well known that if $\theta\in \widehat{\mathbb{O}}_\mathbb{K}$, the discriminant  of $F_\theta(x)$ has the form 
\begin{equation}\label{equ1}
D(\theta)=I(\theta)^2D_\mathbb{K},
\end{equation}

 where $I(\theta)=(\mathbb{O_K}:\mathbb{Z}[\theta])$ is called the index of $\theta$. Let
 
 \begin{equation}\label{equ2}
 I(\mathbb{K})=\underset{\theta\in \widehat{\mathbb{O}}_\mathbb{K}}{\text{gcd}}\;I(\theta).
\end{equation} 
A prime number $p$ is called a common factor of indices   if $p\mid I(\mathbb{K})$. \\

 The knowledge of the index $ I(\mathbb{K})$ makes it possible to find the explicit decomposition of the prime numbers in the fields  $\mathbb K$ : \   If the prime number $p\nmid I(\mathbb{K})$,  by equation \eqref{equ2}  there exist a primitive integer $\theta$ where $p\nmid I(\theta)$ and by Dedekind's theorem \cite[\S 18]{Delone} we  explicitly  have the factorization of $p$ using $\theta$.  However, if $p$ is a common factor divisor, the prime ideal decomposition in $\mathbb{O_K}$ is more difficult. \\
 
Let us recall the statement of Dedekind's theorem.  Let $\mathbb{K}=\mathbb{Q}(\theta)$ be an algebraic number field with $\theta\in\mathbb{O_\mathbb{K}}$. Let $p$ be a rational prime. Let $$f(x)=\textrm{Irr}_\mathbb{Q}(x, \theta)\in\mathbb{Z}[x].$$
We consider  the canonical surjection map $ 
\mathbb{Z}[x]\to \mathbb{Z}/p\mathbb{Z}[x]$.  We write $$\bar{f}(x)=g_1(x)^{e_1}\cdots g_r(x)^{e_r},$$
where $g_1(x),\cdots,g_r(x)$ are distinct monic irreducible polynomials in $\mathbb{Z}/p\mathbb{Z}[x]$ and $e_1,\cdots,e_r$ are positive integers. \\ 
For $i=1,2,\ldots,r$  denote by $f_i(x)$ any monic polynomial of $\mathbb{Z}[x]$ such that $\bar{f}_i=g_i$. We then set  $$P_i=<p,f_i(\theta)>. $$
If $I(\theta)\not\equiv0\mod p$ then we have $P_1,\cdots,P_r$ are distinct prime ideals of \  $\mathbb{O_\mathbb{K}}$ with 
\[ p\mathbb{O}_\mathbb{K}={P_1^{e_1}}\cdots{P_r^{e_r}},\]
 and 
 \[ N(P_i)=p^{\deg(f_i)} . \]
 \\
 In the following we review some known results on the computation of the index $I(\mathbb K)$.\\
 
$\bullet$  If $\mathbb{K}$ is a quadratic field $(n=2)$, by classical number theory on quadratic fields, one can show that $I(\mathbb{K})=1$. \\
$\bullet$ In case $\mathbb{K}$ is a cubic field $(n=3)$ , Engstrom \cite{Engstrom} showed that $I(\mathbb{K})=1$ or $2$.  
Llorente and Nart  in \cite[Theorem 1]{Nart}  determine the {\it type of decomposition} of the rational primes and in   \cite[Theorem 4]{Nart}   give a necessary and sufficient condition for the index of $\mathbb{K}$ to be $2.$ Moreover,  
in the paper \cite{K. Spearman } Spearman and Williams  give the {\it explicit prime ideal factorization} of $2$ in cubic fields with index $2$.\\
$\bullet$  In the case  $\mathbb{K}$ is a quartic field $(n=4)$,  Engstrom \cite{Engstrom} showed that 
\begin{equation}\label{equ3}
   I(\mathbb{K})=1,2,3,4,6\;\text{or}\;12.
\end{equation}   
$\bullet$ If  the field   $\mathbb{K}$ is a cyclic quartic field , Spearman and Williams \cite{Williams } showed that $I(\mathbb{K})$ assumes all of these values and they give necessary and sufficient conditions for each to occur and find  an asymptotic formula for the number of cyclic quartic fields with discriminant $\leq x$ and $I(\mathbb{K})=i$ for each $i\in \{1,2,3,4,6,12\}$,
$$N(x;i)=\alpha_ix^\frac{1}{2}+O(x^\frac{1}{3}log^3x),$$
where
\begin{align*}
\alpha_1&\approx0.0970153,\quad \alpha_2\approx0.0067627,\quad \alpha_3\approx0.0101764,\\
\alpha_4&\approx0.0067627,\quad \alpha_6\approx0.0006321,\quad \alpha_{12}\approx0.0006321.
\end{align*}

$\bullet$ Funakura \cite[Theorem 5]{Funakura} showed that in the case of a pure quartic field, we have $I(\mathbb{K})=1,2$.\\
$\bullet$ In the pure quartic field ,  Spearman and Williams \cite{K. Williams }  gives the explicit prime ideal factorization of $2$ when the index is equal to $2$ . \\
   \subsection{Counting discriminants with given index}
   The general problem in algebraic number theory is that of counting the number of fields by discriminant. If we let $N_n(x) $ denote the number of fields of degree $n$ over $\mathbb Q$ whose discriminants do not exceed $x$ in absolute value. Then there is a conjecture that 
   \[\lim_{x\to\infty}\frac{N_n(x)}{x}\] exist and is non-zero. This conjecture is proved in many special cases:\\
   
  $\bullet$  For quadratic field, the problem is simple we have
   \[N_2(x)=\sum_{|d|\leq x\atop d\equiv 1\;mod\; 4}1+\sum_{|d|\leq \frac{x}4\atop d\equiv 2,3\;mod\; 4}1\sim \frac{6}{\pi ^2} x.
   \]
   $\bullet$ For cyclic cubic fields H. Cohn \cite{Harvey} and  Cohen-Diaz y Diaz  in \cite[p.577, \S3]{cohen} showed that
   \[N_3(x)\sim C_3 x^{1/2}\]
   where 
   \[C_3=  \frac{11\sqrt{3}}{36\pi}\prod_{p\equiv 1\;mod\; 6\atop p\textrm{ prime }}\left(1-\frac{2}{p(p+1)}\right) .\]
   
   $\bullet$ For cyclic quartic fields Baily \cite[p. 209, Theorem 9]{Baily} and its revised form by many authors see \cite{Maki},  \cite{Ou} and \cite[p.580, \S5]{cohen} by differents methods, as follows
   \[N_4(x)\sim C_4 x^{1/2}\]
   where 
  \[C_4=  \frac{3}{\pi^2}
  \left(
  (1+\frac{\sqrt{2}}{24})   \prod_{p\equiv 1\;mod\;4\atop p\textrm{ prime }}  \left(1+\frac{2}{p^{3/2}+p^{1/2}}\right) -1
  \right).
  \]
   
   \subsection{Summary of the paper}
In this paper, we investigate the family of the simplest quartic fields which are defined by adjunction to $\mathbb{Q}$ of a root of  the polynomial   $$P_m=x^4-mx^3-6x^2+mx+1,\quad m\in\mathbb{Z}^+$$
   where $m^2+16$ is not divisible by an odd square.  
It is easy to see that  $m$ may be specified greater than zero as $P_m$ and $P_{-m}$  generate the same extension. \\

 M.N. Gras in \cite{Gras} proved that  those polynomials are reducibles precisely when $m^2+16$ is a square, which occurs only for excluded cases $m=0,\,3$, and show that the form $m^2+16$ represents infinitely many square-free integers.
Olajos \cite{Olajos} proved that $\mathbb{K}_m$ admits power integral bases only for $m=2,\,4$ and he gave all generators of power integral bases. Recently, Ga\'al and Pentr\'ayani \cite{Gaal} compute the  minimal index of the simplest quartic fields.\\

 The purpose of this paper :\\
 $\bullet$ We study the common indices $I(\mathbb{K}_m)$ and the prime ideal decomposition is determined explicitly, \\
 $\bullet$ We establish an asymptotic formula for the number of the simplest quartic fields with given index and the discriminants less than $x$.
\section{\textbf {Statement of main results}}
We state our main results.
\subsection{Computation of the index $I(\mathbb{K}_m)$ and prime decomposition}
\begin{theorem}\label{thm1} Let $m\in\mathbb{Z}^+,\;m\neq 0,3$ and $\mathbb{K}_m=\mathbb{Q}(\theta)$ where $\theta$ be a root of $P_m$. 
Then we have
 \begin{enumerate}
\item 

\[I(\mathbb{K}_m)=
\begin{cases}
   2 & \text{if \; $m$ odd},\\
   1 & \text{if \; $m$ even}.   
\end{cases}
\]

\item The prime ideal factorization of $2$ in $\mathbb{K}_m$ with index $2$ is
$$2\mathbb{O_K}_m=<2,\frac{1+\sqrt{m^2+16}}{2}><2,\frac{1-\sqrt{m^2+16}}{2}>.$$
\item The prime ideal factorization of $2$ in $\mathbb{K}_m$ with index $1$ is as follows:
\begin{enumerate}
\item[1)] If  $v_2(m)=1$, then,
$$2\mathbb{O_K}_m=<2,\frac{m\theta^3+10\theta^2-m\theta+6}{4}>^2.$$
\item[2)] If $v_2(m)=2$, then,
$$2\mathbb{O_K}_m=<2,\frac{1+\theta+\theta^2+\theta^3}{4}>^4.$$
\item[3)] If $v_2(m)=3$, then,
$$2\mathbb{O_K}_m=<2,\frac{(m^3+25m+4)\theta^3+(5m^2+168)\theta^2-(m^3+21m-56)\theta-m^2+8}{16}>^2.$$
\item[4)] If $v_2(m)\geq 4$, then,
$$2\mathbb{O_K}_m=<2,\frac{2+7\theta+\theta^3}{4}>^2<2,\frac{5+7\theta+\theta^3}{4}>^2.$$
\end{enumerate}
\end{enumerate}
\end{theorem}
\begin{proposition}\label{Pro1}
Let $m\in\mathbb{Z}^+,\;m\neq 0,3$ and $\mathbb{K}_m=\mathbb{Q}(\theta)$ where $\theta$ be a root of $P_m$. Then
\begin{center}
$I(\theta)=
\begin{cases}
   2 & \text{if \; $v_2(m)=0$},\\
   2^2 & \text{if \; $v_2(m)=1$ },\\
   2^3 & \text{if \; $v_2(m)=2$},\\
   2^4 & \text{if \; $v_2(m)\geq3$}.\\
\end{cases}$
 \end{center}
 \end{proposition} 
\begin{remark}
Note that we can show that for any $p$ odd prime, we acn 
give explicitly the factorization of $p$ using $\theta$. To do it, we need to factorize $P_m\;mod\;p$ and we use Dedekind's theorem.
\end{remark}
\subsection{Asymptotic number of simplest quartic fields with discriminant $\leq x$ and given index} 
Let $\mathbb{K}_m=\mathbb{Q}(\theta)$ where $\theta$ be a root of $P_m$ (simplest quartic fields). \\
We define for a positive integer $i$
\begin{equation}\label{equ5}
N(x,i)=\text{number of $\mathbb{K}_m$ with $D(\mathbb{K}_m)\leq x$ and $I(\mathbb{K}_m)=i$}.
\end{equation}
We state our second main result.
\begin{theorem}\label{thm2} For $i=1,2,3,4,6,12$ we have
\begin{enumerate}
\item $N(x,i)=0$  for $ i\neq1,2$.
\item For $i=1,2$, we have the asymptotic formulas
\[
N(x,i)\sim  C_ix^{\frac{1}{6}},
\]
where 
\[
C_1=\left(\frac{1}{4\sqrt[3]{4}}+\frac{1}{4\sqrt[3]{2}}+\frac{1}{4}\right) \prod_{p\equiv 1\;mod\;4\atop p\textrm{ prime }} \left(1-\frac{2}{p^2}\right)
\]
and 
\[
C_2=\frac{1}{4}\prod_{p\equiv 1\;mod\;4\atop p\textrm{ prime }} \left(1-\frac{2}{p^2}\right)
\]
\end{enumerate}

\end{theorem}
As corollary, from Theorem \ref{thm2}, the number $N(x)$ of simplest quartic fields $\mathbb{K}_m$ with discriminant $D(\mathbb{K}_m)\leq x$ is given by 
$$N(x)= N(x,1)+N(x,2)\sim \left(\frac{1}{4\sqrt[3]{4}}+\frac{1}{4\sqrt[3]{2}}+\frac{1}{2}\right) \prod_{p\equiv 1\;mod\;4\atop p\textrm{ prime }} \left(1-\frac{2}{p^2}\right)\ x^{\frac{1}{6}} , \textrm{ as } x\to \infty$$

\end{document}